\documentclass[11pt]{article}
\usepackage{amssymb,amsmath,amsthm}
\newtheorem{theorem}{Theorem}

\newtheorem{proposition}{Proposition}

\date{}
\numberwithin{equation}{section} \numberwithin{theorem}{section}
\numberwithin{lemma}{section} \numberwithin{corollary}{section}
\numberwithin{remark}{section} \numberwithin{proposition}{section}
\numberwithin{definition}{section}
\usepackage{color}
\usepackage{curves}
\usepackage[dvips]{graphics}
\usepackage{amsfonts}
\usepackage{amsmath}
\usepackage{amssymb}
\usepackage[normalem]{ulem}
\usepackage{t1enc}
\usepackage{fancybox}
\usepackage{epsfig}
\usepackage{minitoc}
\usepackage{fancyhdr}

\begin{document}

\newcommand{\n}{\noindent}

\newcommand{\vs}{\vskip}

\title{Uniqueness of solution of a heterogeneous evolution dam problem associated with a compressible fluid flow through  a rectangular porous medium  }

\vs 0.5cm
\author{E. Zaouche \\
Department of Mathematics\\
University of EL Oued B. P. 789 El Oued 39000 Algeria\\
E-mail: elmehdi-zaouche@univ-eloued.dz} \maketitle

\begin{abstract}
This paper is concerned with  an uniqueness of solution of the weak
formulation  of an
 evolution dam problem related to a compressible fluid
flow through a two-dimensional, rectangular and heterogeneous porous
medium. Note that our problem associated with the equation
$a(x_1)(u_{x_2}+\chi)_{x_2}-(u+\chi)_t=0$. Our technique is based on
the idea that we transform  the weak form of this equation into  a
similar situation to the proof of the uniqueness in the
incompressible case (see [12]). It is also difficult to adapt the
proof obtained in [12] by using some properties of the solutions as
in [12, Sect. 2].
\end{abstract}

\begin{flushleft}
2010 Mathematics Subject Classification: 35A02, 35B35,  76S05.
\end{flushleft}

\begin{flushleft}
Key words: Heterogeneous evolution dam problem; compressible fluid
flow;  rectangular porous medium; uniqueness.
\end{flushleft}

\section{Introduction }\label{s1}
In this paper,  we consider the following  weak formulation of the
evolution dam problem with heterogeneous coefficients which related
to a compressible fluid flow:
\begin{equation*}{\bf (P)} \quad \left\{
\begin{aligned}
&\text{Find }
(u, \chi) \in L^2(0,T;H^1(\Omega))\times L^\infty (Q) \text{ such that}:\\[0.2cm]
&\quad u \geq 0, \; 0\leq \chi\leq 1,\; u.(1-\chi) = 0  \quad \text{ a.e. in }Q \\[0.1cm]
&\quad u=\phi \hspace{0.34cm} \text{ on } \Sigma_2 \\
&\quad \displaystyle{\int_Q \big[a(x_1)(u_{x_2} +\chi)\xi_{x_2}-( u+\chi)\xi_t\big] dx\,dt}\\
&\qquad\qquad  \leq \displaystyle{\int_\Omega (\chi_0(x)+ u_0(x))\xi(x,0)\,dx}  \\
&\quad \forall \xi \in H^1(Q),\; \xi=0 \text{ on }
   \Sigma_3,\; \xi \geq 0 \text{ on } \Sigma_4,\; \xi(x,T)=0 \,\text{ for a.e. }
x\in \Omega,
\end{aligned} \right.
\end{equation*}
where $\Omega=(0,L)\times(0,H)$ is a bounded rectangular domain in
$\mathbb{R}^{2}$ that the generic point in $\Omega$ is denoted by
$x=(x_1, x_2)$,  which represents a porous medium,  with  a boundary
$\partial\Omega=\Gamma_1\cup\Gamma_2$ such that $\Gamma_1=
[0,L]\times\{0\}, \, \Gamma_2=(\{0\}\times[0,H])\cup
([0,L]\times\{H\})\cup (\{L\}\times[0,H]).$
 $Q=\Omega\times(0,T), \,T>0$ is a positive real number, $\phi$ is
a nonnegative Lipschitz continuous function defined in
$\overline{Q},\, \Sigma_1=\Gamma_1\times(0,T),\,
\Sigma_2=\Gamma_2\times(0,T),\, \Sigma_3=\Sigma_2\cap \{\phi>0\}$
and $\Sigma_4=\Sigma_2\cap \{\phi=0\}$. Note that from a physical
point of view, $\Gamma_1$ is  the impervious part of
$\partial\Omega$,  $\Gamma_2$ is the part in contact with either air
or the fluid reservoirs and $\phi$ represents the assigned pressure
on $\Sigma_2.$
 \n For some two constants $0<\lambda\leq \Lambda, \, a(x_1)$ is a
 function of the variable $x_1$ satisfying
\begin{eqnarray}\label{e1.3-1.4}
&&\lambda \leq a(x_1)\leq \Lambda\quad\text{ a.e.} \;
x_1\in(0,L)\nonumber
\end{eqnarray}
and $u_0, \chi_{0}: \Omega\rightarrow \mathbb{R}$ such that for some
positive constant $M,$
\begin{eqnarray}\label{e1.3-1.4}
&&0\leq\chi_{0}(x)\leq 1 \quad\text{ a.e.} \; x\in\Omega,\nonumber\\
&&0\leq u_{0}(x)\leq M \quad\text{ a.e.} \; x\in\Omega.\nonumber
\end{eqnarray}
Note that the strong formulation of $(P)$ associated with the
initial data $(u_0,\chi_0)$ is  given by
\begin{equation*}{\bf (SF)} \quad \left\{
\begin{aligned}
\quad u \geq 0, \; 0\leq \chi\leq 1,\; u(1-\chi) &= 0
&\quad & \text{ in }Q \\
a(x_1)(u_{x_2} +\chi )_{x_2}-(u+\chi)_t &= 0 & \quad & \text{ in } Q \\
u &=\phi&\quad & \text{ on } \Sigma_2 \\
u(\cdot,0)+\chi(\cdot,0)&=u_0+\chi_0 & \quad & \text{ in } \Omega \\
a(x_1)(u_{x_2} +\chi )\cdot \nu &=0  & \quad
& \text{ on } \Sigma_1 \\
a(x_1)(u_{x_2} +\chi )\cdot \nu &\leq 0  & \quad & \text{ on }
\Sigma_4.
\end{aligned} \right. \hspace{3.3cm}
\end{equation*}

\n Concerning existence of a solution of  problem $(P)$, we refer to
[2], where an existence theorem was established for a class of
non-stationary  free boundary problems including the problem $(P)$.
Also, a regularity in time for this class was obtained in [9].

\n Uniqueness of the solution was proved in [3] and [12] by using
the method of doubling variables respectively for a homogeneous dam
with general geometry and for rectangular and heterogeneous porous
medium related to an incompressible fluid flow. By a different
method, uniqueness was obtained in [5] and [10] for a rectangular
dam wet at the bottom and dry near to the top, respectively, in
homogeneous and heterogeneous cases, in both compressible and
incompressible fluids. Also, we refer to [1], for the uniqueness of
the solution of the evolution free boundary problem in theory of
lubrication.

\n In this article, we prove the  uniqueness of the solution of
problem $(P)$ by transforming the weak form of the equation
$a(x_1)(u_{x_2}+\chi)_{x_2}-(u+\chi)_t=0$ in $Q$ into a form does
not contain the temporal term in $u$ (see Proposition 2.2) which
allows us to apply the proof  of the uniqueness in the
incompressible case (see [12]). Note that we can establish some
properties of the solutions as in [12, Sect. 2], but, it is
difficult to adapt the proof of the uniqueness obtained in [12] with
these properties.

\section{Uniqueness of the solution}\label{2}
We assume throughout this section that
$$a\in C^{1}([0,L]).$$

\n First, we have the following regularity result for solutions of
problem $(P)$ ([9]).

\begin{proposition}\label{prop2.2}
We have
\begin{eqnarray}\label{e2.9}
&&\chi\in C^0([0,T];L^{p}(\Omega)), \quad \forall p\in [1,+\infty),\nonumber\\
&&u\in C^0([0,T];L^{p}(\Omega)), \quad \forall p\in [1,2].\nonumber
\end{eqnarray}
\end{proposition}

\n The following proposition plays an important role in the proof of
the uniqueness of  the solution.
\begin{proposition}\label{prop2.2}
Let $(u, \chi)$  be a solution of problem $(P)$. Then we have
\begin{eqnarray}
\forall \xi \in H^{1}_{0}(\Omega), \,\forall t\in [0,T]: \quad
\int_{\Omega}a(x_1) (u_{x_2}+\chi)\xi_{x_2}\,dx=0.
\end{eqnarray}
\end{proposition}
\n \emph{Proof.} Without loss of generality, we assume that $\xi\in
\mathcal{D}(\Omega)$. Then we obtain the result by approximation for
$\xi\in H_{0}^{1}(\Omega).$ Let $\kappa$ be a fixed element of
$(0,T]$. We define the following function $\eta$ on $[0,s]:$
\[ \eta(t)= \begin{cases}
2(\frac{t}{\delta})^2  &\text{if } t\in [0,\frac{\delta}{2}]\\
1-2(1-\frac{t}{\delta})^2 &\text{ if } t\in (\frac{\delta}{2},\delta]\\
1 &\text{if }  t\in (\delta,\kappa-\delta]  \\
1-2(1-\frac{\kappa-t}{\delta})^2 &\text{if } t\in (\kappa-\delta,\kappa-\frac{\delta}{2}]\\
2(\frac{\kappa-t}{\delta})^2  &\text{if } t\in
(\kappa-\frac{\delta}{2},\kappa].
 \end{cases}
\]
where $\delta>0$  is a positive real number.  We see that $\eta \in
C^{1}([0,\kappa]), \, \eta(0)=\eta(\kappa)=0$ and
\[
\eta'(t)= \begin{cases}
4\frac{t}{\delta^2}  &\text{if } t\in [0,\frac{\delta}{2}]\\
\frac{4}{\delta}(1-\frac{t}{\delta}) &\text{if } t\in (\frac{\delta}{2},\delta]\\
0 &\text{if }  t\in (\delta,\kappa-\delta]  \\
-\frac{4}{\delta}(1-\frac{\kappa-t}{\delta})
&\text{ if } t\in (\kappa-\delta,\kappa-\frac{\delta}{2}]\\
-\frac{4}{\delta}(\frac{\kappa-t}{\delta}) &\text{if } t\in
(\kappa-\frac{\delta}{2},\kappa].
 \end{cases}
\]
For all $\xi\in \mathcal{D}(\Omega)$, the function $\xi\eta^2$
belongs to $H^{1}(Q)$ and satisfies $\xi\eta^2=0$ on
$\partial\Omega\times(0,T)$ and
$(\xi\eta^2)(.,0)=(\xi\eta^2)(.,T)=0$ a.e. in $\Omega.$ Therefore,
$\pm \xi\eta^2$ are test functions for $(P)$ and we have
\begin{eqnarray}
&0=&\int_{\Omega\times(0,\kappa)}[a(x_1)
(u_{x_2}+\chi)\xi_{x_2}\eta^2\,dx dt
-\int_{\Omega\times(0,\kappa)}2(u+\chi)\eta\eta'\xi\,dx
dt=0\label{e3.44}\nonumber\\
&:=&I_{\delta}^{1}-I_{\delta}^{2}.
\end{eqnarray}
Applying the Lebesgue dominated convergence theorem to
$I_{\delta}^{1}$, we obtain
\begin{equation}
\lim_{\delta\to0}I_\delta^{1}
=\int_{\Omega\times(0,\kappa)}a(x_1)(u_{x_2}+\chi)\xi_{x_2}\,dx
dt.\label{e3.46}
\end{equation}
Moreover, we use the definition of $\eta'$, we see that the quantity
$|I_{\delta}^{2}|$ can be estimated as follows:
\begin{equation}\label{e3.47}
\begin{aligned}
|I_{\delta}^2|
=&2\Big|\int_{\Omega}\int_{0}^{\kappa}a(x_1)(u+\chi)\eta\eta'\xi
\,dx dt
+\int_{\Omega}\int_{\kappa-\delta}^{\kappa}a(x_1)(u+\chi)\eta\eta'\xi
\,dx dt\Big| \\
\leq& C\Big\{ \int_{0}^{\delta}|u+\chi|_{1,\Omega}\eta|\eta'|\, dt
+\int_{\kappa-\delta}^{\kappa}|u+\chi|_{1,\Omega}\eta|\eta'|  \,dt\Big\} \\
:=&C(I_{\delta}^{2,1}+I_{\delta}^{2,2})
\end{aligned}
\end{equation}
where $C=\sup_{(x_1,x_2)\in \Omega}(a(x_1)|\xi(x_1,x_2)|)$. Let us
set $f(t)=|u+\chi|_{1,\Omega}\eta.$ We see that $f$ is
right-continuous and vanishes at $0$ since $u, \chi \in C^{0}([0,T];
L^{1}(\Omega)) $\,(see Proposition 2.1), \,$ \eta\in
C^{0}([0,\kappa])$ and $\eta(0)=0.$  So, by using $|\eta'|\sim
\frac{1}{\delta}$, we get
\begin{equation}\label{e3.48}
\lim_{\delta\to0}I_{\delta}^{2,1}=0.
\end{equation}
Similarly, the function $f$ is  left-continuous and vanishes at
$\kappa$ and $|\eta'|\sim \frac{1}{\delta}$ to obtain
\begin{equation}\label{e3.48}
\lim_{\delta\to0}I_{\delta}^{2,2}=0.
\end{equation}
Hence, by letting $\delta\rightarrow 0$ in (2.4) and using
(2.5)-(2.6), we get
\begin{equation}\label{e3.48}
\lim_{\delta\to0}I_{\delta}^{2}=0.
\end{equation}
Now, using (2.3) and (2.7) we obtain by taking $\delta\rightarrow 0$
in (2.2):
\begin{eqnarray}
\forall \kappa\in (0,T]: \quad\int_{\Omega\times(0,\kappa)}a(x_1)
(u_{x_2}+\chi)\xi_{x_2}\,dx dt =0\nonumber
\end{eqnarray}
which leads, using integration by parts, to
\begin{eqnarray}
\forall \kappa\in (0,T]: \quad0=\int_{0}^{\kappa}\int_{\Omega}a(x_1)
(\chi\xi_{x_2}-u\xi_{x_2x_2})\,dx dt :=F(\kappa).\nonumber
\end{eqnarray}
Then, we have
\begin{eqnarray}
\forall t\in [0,T]: \quad F^{\prime}(t)=\int_{\Omega}a(x_1)
(\chi\xi_{x_2}-u\xi_{x_2x_2})\,dx=0 \nonumber
\end{eqnarray}
since $t\mapsto \displaystyle{\int_{\Omega}a(x_1)
(\chi\xi_{x_2}-u\xi_{x_2x_2})\,dx} $ is continuous on $[0,T].$ Using
again integration by parts, we obtain (2.1).\qed

\n Let us now give the uniqueness theorem.
\begin{theorem}\label{thm3.1}
Let $(u_{1}, \chi_{1})$ and $(u_{2}, \chi_{2})$ be two solutions of
problem $(P)$ such that
$u_1(.,0)+\chi_1(.,0)=u_2(.,0)+\chi_2(.,0)=u_0+\chi_0$ a.e. in
$\Omega$. Then we have
$$(u_1,\chi_1)=(u_2,\chi_2) \quad \text{a.e.
in } Q.
$$
\end{theorem}
\n \emph{Proof.} Let $\xi\in \mathcal{D}(\Omega)$ and $\eta\in
\mathcal{D}(0,T)$ such that $\xi\geq0$ and $\eta\geq0.$ For a
positive real number $\delta>0$, we consider the following functions
 $\rho_{1,\delta}(r)=\frac{1}{\delta}\rho_{1}(\frac{r}{\delta})$,
$\rho_{2,\delta}(r) =\frac{1}{\delta}\rho_{2}(\frac{r}{\delta})$,
$\rho_{3,\delta}(r) =\frac{1}{\delta}\rho_{3}(\frac{r}{\delta}) $
 with $\rho_{1}, \rho_{2}, \rho_{3}\in \mathcal{D}(\mathbb{R})$,
$ \rho_{1}, \,\rho_{2}, \,\rho_{3}\geq0$,
$\operatorname{supp}(\rho_{1}), \,\operatorname{supp}(\rho_{2}),\,
\operatorname{supp}(\rho_{3})\subset(-1,1)$. We define
\begin{eqnarray}\label{e2.9}
&&\forall (x,t,y,s)\in \overline{Q\times Q}: \nonumber\\
&&\zeta(x,t,y,s)=\xi(\frac{x_{1}+y_{1}}{2},\frac{x_{2}+y_{2}}{2})\eta(\frac{t+s}{2})
\rho_{1,\delta}(\frac{x_{1}-y_{1}}{2})\rho_{2,\delta}(\frac{x_{2}-y_{2}}{2})
\rho_{3,\delta}(\frac{t-s}{2}).\nonumber
\end{eqnarray}
 By choosing $\delta$ small enough, we obtain
\begin{gather}
\forall (t,y,s)\in (0,T)\times Q: \quad\zeta(\cdot,t, y,s)=0  \quad \text{ on } \partial \Omega,\label{e3.3}\\
\forall (x,t,s)\in Q\times(0,T): \quad \zeta(x,t, \cdot,s) \quad
\text{ on }
\partial \Omega.\label{e3.4}
\end{gather}
With this notation of doubling variables, we consider the solution
$(u_1,\chi_1)$ is in the variable $(x,t)$ and $(u_2,\chi_2)$  is so
related to the variable $(y,s).$  So, for a positive real number
$\epsilon>0$, we set
\begin{equation}
\vartheta(x,t,y,s)=\min\Big({{(u_{1}(x,t)-u_{2}(y,s))^+}\over
\epsilon},\zeta(x,t,y,s)\Big).\label{e3.5}
\end{equation}
From (2.8), the function $\xi(.)=\vartheta(.,t,y,s)\in
H_{0}^{1}(\Omega)$ for all  $(t,y,s)\in(0,T)\times Q$. Therefore, we
can apply Proposition 2.2 to $(u_1,\chi_1)$ with
$\xi(.)=\vartheta(.,t,y,s)$ to get
\begin{equation}\label{e3.7}
\int_{\Omega} a(x_1)(u_{1x_2}+\chi_1 )\vartheta_{x_2}\,dx=0
\end{equation}
 and
by integrating
 over $(0,T)\times Q$, we obtain
\begin{equation}\label{e3.7}
\int_{Q\times Q} a(x_1)(u_{1x_2}+\chi_1 )\vartheta_{x_2}\,dx dt dy
ds=0.\nonumber
\end{equation}
On the other hand,  if we use $u_1(1-\chi_1)=0$ a.e. in $Q$, one can
easily verify that for a.e. $(x,t,y,s)\in Q\times Q:$
\begin{equation}
\chi_{1} a(x_1)\Big(\min\Big({{(u_{1}-u_{2})^+}\over
\epsilon},\zeta\Big)\Big)_{x_2}=a(x_1)\Big(\min\Big({{(u_{1}-u_{2})^+}\over
\epsilon},\zeta\Big)\Big)_{x_2}.\nonumber
\end{equation}
Therefore,  (2.11) takes on the form
\begin{equation}\label{e3.7}
\int_{Q\times Q} a(x_1)(u_{1x_2}+1)\vartheta_{x_2}\,dx dt dy ds=0.
\end{equation}
Similarly, by (2.9), we see that for all $(x,t,s)\in Q\times(0,T)$,
we can apply Proposition 2.2 to  $(u_2,\chi_2)$ with
$\xi(.)=\vartheta(x,t,.,s)$ to obtain
\begin{equation}\label{e3.7}
\int_{\Omega}a(y_1)(u_{2y_2}+\chi_2)\vartheta_{y_2}\, dy=0\nonumber
\end{equation}
 and by integrating over $Q\times(0,T)$, we have
\begin{equation}\label{e3.7}
\int_{Q\times Q}a(y_1)(u_{2y_2}+\chi_2)\vartheta_{y_2}\, dx dt dy
ds=0.
\end{equation}
By subtracting (2.13) from (2.12), we find
\begin{equation}\label{e3.9}
\begin{aligned}
&\int_{Q\times Q} \big[a(x_1) u_{1x_2}\vartheta_{x_2}-a(y_1)
u_{2y_2}\vartheta_{y_2} \\
&+a(x_1)\vartheta_{x_2}-\chi_{2}a(y_1)\vartheta_{y_2}\big]\,dx dt dy
ds=0.
\end{aligned}
\end{equation}
Now, we see that the relation (2.14) is in a similar situation  to
[[12],  (3.9)]. So, by using the method of doubling variables and
argue as in [12], we arrive at $u_1=u_2$ a.e. in $Q$ (see  the first
part of the proof of [[12], Theorem 3.4]). In order to get
$\chi_1=\chi_2$ a.e. in $Q$, we use an argument similar to the
second part of the proof of Theorem 3.4 of [12].

\end{document}